\documentclass[12pt]{amsart}
\usepackage{geometry}                
\usepackage{graphicx}

\usepackage{amsmath}
\usepackage{amsfonts}
\usepackage{enumerate}
\usepackage{amssymb}
\usepackage{epstopdf}

\theoremstyle{definition}

\numberwithin{equation}{section}

\newtheorem{Thm}{Theorem}[section]

\newtheorem{Pro}[Thm]{Proposition}

\renewcommand{\epsilon}{\varepsilon}

\renewcommand{\rho}{\varrho}

\begin{document}

\title{Balanced presentaions of the trivial group and four-dimensional geometry}

\author{ Boris Lishak and Alexander Nabutovsky }
\maketitle

\begin{abstract}
We prove that 1) There exist infinitely many non-trivial
codimension one ``thick" knots in $\mathbb{R}^5$; 2) For each
closed four-dimensional smooth manifold $M$ and 
for each sufficiently small positive $\epsilon$ the set of isometry classes
of Riemannian metrics with volume equal to $1$ and injectivity radius greater than $\epsilon$ is disconnected; 3) For each closed four-dimensional
$PL$-manifold $M$ and any $m$ there exist arbitrarily large values of $N$ such that
some two triangulations of $M$ with $<N$ simplices cannot be connected
by any sequence of $<M_m(N)$ bistellar transformations, where
$M_m(N)=\exp(\exp(\ldots \exp (N)))$ ($m$ times).
\end{abstract}

\section{Main results.} \label{Main}

{\bf 1.1.} The goal of this paper is to extend results of [N1], [N2], [N3]
to the four-dimensional situation.

\begin{Thm} Let $M$ be any closed four-dimensional Riemannian manifold.
Let $I_{\epsilon}(M)$ denote the space of isometry classes of Riemannian
metrics on $M$ with volume equal to $1$ and injectivity radius greater than
$\epsilon$. (This space is endowed with the Gromov-Hausdorff metric $d_{GH}$.)
Then for all sufficiently small $\epsilon >0$ $I_\epsilon(M)$ is disconnected,
and, moreover, can be represented at the union of two non-empty subsets
$A_1,\ A_2$ such that for any $\mu_1\in A_1$, $\mu_2\in A_2$
$d_{GH}(\mu_1,\ \mu_2)>{\epsilon\over 10}$.
\par
Furthermore, let for each $m$ $\exp_m(x)$ denote $\exp(\exp(\ldots(\exp x)))$ ($m$ times).
Then for each $m$ for all sufficiently small $\epsilon$ there exist $\mu,\nu\in I_{\epsilon}(M)$ with the following property. Let
$\mu_1=\mu, \mu_2,\ldots , \mu_N=\nu$ be a sequence of isometry classes of Riemannian metrics on $M$
of volume one such that for each $i$ $d_{GH}(\mu_i,\ \mu_{i+1})\leq {\epsilon\over 10}$. Then
$\inf_i inj(\mu_i)\leq {1\over \exp_m({1\over\epsilon})}$.
\end{Thm}

A (stronger) analog of this theorem for $n>4$ as well as for a class of closed four-dimensional manifolds
representable as the connected sum of any closed $4$-manifold and several copies of $S^2\times S^2$ can be found in [N2] (Theorem 1 and section 5.A).
(More precisely, ``several" means $14$. The minimal number of copies of $S^2\times S^2$ required for the method
of [N2] to work is equal to the number of relators in a sequence of finitely presented groups, where the triviality problem is algorithmically unsolvable;
cf. [Sh1], [Sh2]).

The next theorem is a four-dimensional analog of Theorem 11 from [N2].
For each smooth manifold $M$ define $Al_1(M)$ as the space of $C^1$-smooth
Alexandrov spaces of curvature $-1\leq K\leq 1$, $C^1$-diffeomorphic to $M$
(cf. [BN] for a definition of Alexandrov spaces with two-sided bounds on
sectional curvature).
A result of I. Nikolaev ([Ni]) implies that all of them are Gromov-Hausdorff
limits of sequences of smooth Riemannian structures on $M$. The classical Gromov-Cheeger compactness theorem implies that all elements of $Al_1(M)$ are
$C^{1,\alpha}$-smooth Riemannian structures on $M$ for any $\alpha\in (0,1)$.
We can consider diameter as a functional on $Al_1(M)$.

\begin{Thm} Let $M$ be a closed $4$-dimensional manifold such that 
either its Euler characteristic is not equal to zero, or its simplicial
volume is not equal to zero. Then diameter regarded as a functional on $Al_1(M)$has infinitely many local minima. The set of values of diameter at its local minima on $Al_1(M)$ is unbounded. 
\end{Thm}

The assumptions about $M$ imply a uniform positive lower bound for the volume of
all elements of $Al_1(M)$. Now the Gromov-Cheeger theorem implies
the compactness of sublevel sets of $diam$, $diam^{-1}((0,x])$, on $Al_1(M)$
for all values of $x$. Now we see that it is sufficient to prove that there
exists an unbounded sequence of values of $x$ such that the set of all smooth
Riemannian structures on $M$ with $-1\leq K\leq 1$ and $diam\leq x$
is disconnected, and, moreover, can be represented as a union of two non-empty
subsets with disjoint closures. After noticing that the classical
Cheeger inequality implies that for all such smooth Riemannian structures
the injectivity radius will be bounded below by an explicit positive
function of $x$ (that behaves as $const\ \exp(-3x)$) we see that this
theorem is similar to the previous one, and, in fact, has a very similar
proof.

Theorem 11 in [N2] should not be confused with a much deeper and significantly
more difficult main theorem in [NW1] (see also [NW2] and [W]) that
does not have the assumption that a smooth manifold $M$ of dimension
greater than four has either a non-zero Euler characteristic
or a non-zero simplicial volume, and, therefore, one lacks an a priori uniform
positive lower bound for the volumes of the considered metrics.
At the moment we are
not able to prove a four-dimensional analog of the main theorem of [NW1].

{\bf 1.2.} In order to state the next theorem define {\it crumpledness} 
(a.k.a {\it ropelength}) of an embedded
closed manifold $X^n$ in a complete Riemannian manifold $Y^{n+k}$
as $\kappa(X^n)={vol^{1\over n}(X^n)\over r(X^n)}$, where $r(X^n)$ denotes the
injectivity
radius of the normal exponential map of $X^n$. Informally speaking, $r(X^n)$
can be interpreted as the smallest radius of a nonself-intersecting tube
around $X^n$. This functional was defined in [N1] for hypersurfaces
and named ``crumpledness", but  in later papers on ``thick" knots in $\mathbb{R}^3$ it had been given a new name ``ropelength", as it can be interpreted
as the length of a similar knot such that the maximal radius
of a nonself-intersecting tube around this knot is equal to one (i.e.
it is the length of a similar knot tied on ``thick" rope of radius one).
One of the ideas of [N1] was that one can similarly consider higher-dimensional
``thick" knots. Two knots (=embeddings of $S^n$ in $\mathbb{R}^{n+k}$)
belong to the same $x$-thick knot type
if they both are in the same path component of the sublevel set 
$\kappa^{-1}((0,x])$ of  $\kappa$. 
\par
To state our main result about ``thick" knots it is convenient
to first introduce a space of non-parametrized $C^{1,1}$-smooth embeddings 
$E_n=Emb(S^n, R^{n+1})/Diff(S^n)$ of $S^n$ into $\mathbb{R}^{n+1}$,
and then define $Knot_{n,1}$ as the quotient of $E_n$ with respect
to the action of the group generated by isometries and homotheties
of the ambient Euclidean space $\mathbb{R}^{n+1}$.
The choice of smoothness is motivated by the facts that 1) $r(\Sigma^n)>0$
for every $C^{1,1}$-smooth closed hypersurface; 2) $r$ is an upper semi-continuous functional of $E_n$ and, therefore, $Knot_{n,1}$ (Theorem 5.1.1 of
[N1]); and 3) The sublevel sets
$\kappa^{-1}((0,x])$ in $Knot_{n,1}$ are compact (see [N1], proof of
Theorem 5.2.1). These facts are true for all dimensions $n$.
Now for each $x$ we can consider ``thick" knot $x$-types as subsets of either $E_n$ or $Knot_{n,1}$. A knot $x_1$-type
and $x_2$-type are {\it distinct} if they do not intersect in $E_n$ (or, equivalently, in $Knot_{n,1}$). (Assuming that, say, $x_1\leq x_2$, this is 
equivalent to the $x_1$-knot not being a subset of the $x_2$-knot.)
Our next results imply that there exist non-trivial types of ``thick" four-dimensional knot types of codimension one.

\begin{Thm} 
There exists an infinite sequence of distinct $x_i$-knot types in $E_4$ (correspondingly, $Knot_{4,1}$),
where $x_i$ is an unbounded increasing sequence. Moreover, there exists an unbounded increasing sequence
of $x_i$, which are the values of $\kappa$ at its local minima $k_i$ on $E_4$ (or, equivalently,  $Knot_{4,1}$).
Further, for each $m$ one can find such a sequence of numbers $\{x_i\}$ and knots $k_i$ with the additional
property that any isotopy between $k_i$ and the standard $4$-sphere of radius $1$ in $\mathbb{R}^5$ must pass through
hypersurfaces, where the value of $\kappa$ is greater than $\exp_m(x_i)$. 
\end{Thm}

\noindent
{\bf Remarks.} {\bf 1.} The second assertion of the theorem is stronger than the first assertion, as each local minimum of $\kappa$
with value $x$ gives rise to a $x$-knot type that consists of one knot, if the local minimum is strict, and a connected set
of knots in $\kappa^{-1}(\{x\})$ otherwise. On the other hand, the second asserton immediately follows from the first assertion 
and the compactness of sublevel sets of $\kappa$ (see Theorem 5.1.1 in [N1]).
\par\noindent
{\bf 2.} The local minima of $\kappa$ were called self-clenching hypersurfaces
in [N1]. The idea behind this metaphor is that one can imagine that this
hypersurface is made of very thin material that bends but cannot be stretched.
If it also cannot be squeezed, then the ``thick" hypersurface
is tightly
folded in $\mathbb{R}^5$. It can move (other than rigid body movement) only
if the local minimum is not strict, and only by ``sliding movements", so that
at each moment of time it is still a local minimum of $\kappa$ (i.e. it cannot
be unfolded into a less crumpled shape).
\par\noindent
{\bf 3.} Two very interesting question are whether or not
there exist non-trivial ``thick" knots of codimension one in $\mathbb{R}^3$
and  $\mathbb{R}^4$. The second of these questions is related to 
the smooth Schoenflies conjecture, that asserts that each smooth embedding
of $S^3$ into $\mathbb{R}^4$ is isotopic to the standard round sphere of radius one. (This fact is known for all other dimensions). Note, that if
the smooth Schoenflies conjecture turns out to be false one can still ask
whether or not there are non-trivial ``thick" knot types in the component
of $Knot_{3,1}$ that consists of $3$-spheres in $\mathbb{R}^4$ that are isotopic to the round sphere.  It seems almost ``self-evident" that there are no
non-trivial ``thick" knots $S^1\subset \mathbb{R}^2$, but we do not know a proof of this fact
and are not aware of any publications in this direction.
\par
{\bf 1.3.} To state our third result for every closed four-dimensional $PL$-manifold $M$
consider the set of all simplicial isomorphism classes of simplicial
complexes PL-homeomorphic to $M$. For brevity, we call them {\it triangulations} of $M$. The discrete set $T(M)$ of all triangulations of $M$ can be turned into a metric space using {\it bistellar transformations}. Bistellar transformation
are operations that transform one triangulation into the other as follows.
Let $T_1$ be a triangulation of $M$. Assume that it contains a 
simplicial subcomplex $K$ that
consists of $k$, $1\leq k\leq 5$, $4$-dimensional simplices (together with their faces)
and is simplicially isomorphic to a subcomplex $C$ of the boundary of a
$5$-dimensional simplex $\partial \Delta^5$. To perform the corresponding
bistellar transformation one first removes these $k$ simplices (and all
their faces) and then attaches the closure
of the complement $\partial\Delta^5\setminus C$ to the boundary of $K$
(which is simplicially isomorphic to the boundary of
$\partial\Delta^5\setminus C$).
Since we exchange one PL-disc (triangulated with $k$ $4$-simplices) for 
another (triangulated with $6-k$ $4$-simpleces), we obtain a triangulation
$T_2$ of the same manifold. Moreover, endow
$T_1$ and $T_2$ with length metrics such that each simplex is
a flat regular simplex with side length one. In this case it is easy to see
that
$T_1$ and $T_2$ will be bi-Lipschitz homeomorphic, and the Lipschtz
constants of the homeomorphism and its inverse
will not exceed an absolute constant
that can be
explicitly evaluated. U. Pachner proved that every two triangulations
of the same closed $PL$-manifold can be connected by a finite sequence
of bistellar transformations ([P]). Now one can define the distance
$d_{Bist}(T_1,T_2)$ on $t(M)$ as the minimal number of bistellar
transformations required to transform $T_1$ into $T_2$.

\begin{Thm} For each $4$-dimensional closed PL-manifold $M$ and each
positive integer value of $m$ there
exist arbitarily large values of $N$ and two triangulations $T_1$, $T_2$
with $\leq N$ simplices such that $d_{Bist}(T_1,\ T_2) >\exp_m(N)$.
(In other words, $T_1$ and $T_2$ cannot be connected by any sequence
of less than $\exp_m(N)$ bistellar transformations).
\end{Thm}

A stronger version of this theorem had been proven in [N3] for all manifolds
of dimension greater than four as well as all four-dimensional manifolds that can be represented as
a connected sum with $k$ copies of $S^2\times S^2$, where the value of
$k$ can be chosen as $14$ using
[Sh1], [Sh2].  Note that results of [N3] and, especially, Theorem 1.4
for $M=S^4$
have potential implications for
four-dimensional Euclidean Quantum Gravity (see [N4] and references there).

\section {Proofs.} \label{Proofs}

{\bf 2.1. Balanced finite presentations of the trivial group.} 
In [L] one of the authors have constructed a sequence of finite
balanced presentations of the trivial group. These finite presentations
have two generators and and two relations. They can be described
as the Baumslag-Gersten group $B=<x, t\vert x^{x^t}=x^2>$
with an added variable second relation.
(Here $a^b$ denotes $aba^{-1}$.)
To describe this relation note that there exists a sequence of
words $v_n$ in $B$ of length $O(n)$ representing $x^{E([\log_2\ n])}$,
where $E(m)$ denotes $2^{2^{\dots^2}}$ ($m$ times). Clearly, $v_n$
commutes with all powers of $x$. Words $w_n=[v_n,x]$ represent
the trivial element but one needs to apply the only
relation at least $E([\log_2 n]-const)$ times to establish this fact. (Thus, the Dehn function
of the one relator group $B$ grows faster than any tower of exponentials
of a fixed height of $n$, cf. [Ge], [Pl].)
The extra relation added to $B$ is
$[v_n, x^3][v_n, x^5][v_n, x^7]=t$, where $[a,b]$ denotes the commutator
$aba^{-1}b^{-1}$. The most important property of this sequence of groups
is that any representation of either $x$ or $t$ as a product of conjugates
of the two relators and their inverses will require at least $E([\log_2 n]-2)$
multipliers. Also, note that these finite presentations satisfy the Andrews-Curtis
conjecture. The importance of the last observation is in the fact that
when one constructs a representation complex $K$ of such a finite presentation
(that is, a $2$-complex 
with one $0$-dimensional cell, two $1$-dimensional 
cells corresponding to the generators and two $2$-dimensional cells 
corresponding to the relators), embeds it in $\mathbb{R}^5$,
takes the boundary of a small neighborhood of the embedding, and smoothes
it out, one obtains not merely a smooth homotopy $4$-sphere that must be homeomorphic to $S^4$ by virtue of Freedman's proof of the $4$-dimensional
Poincare conjecture, but a manifold that is diffeomorphic to $S^4$. This fact can be demonstrated without the $4$-dimensional Poincare conjecture using instead
the fact that the operations in the Andrews-Curtis conjecture 
correspond to certain diffeomorphisms of the underlying
manifold (``handle slidings"). A sequence of these diffeomorphisms
corresponding to handle slides will eventually result in the standard sphere
that corresponds to the representation $2$-complex of the trivial
finite presentation of the trivial group (cf. [BHP]).

The resulting smooth hyperspheres in $\mathbb{R}^5$  that
will be denoted by $S^4(v_n)$ can, after rescaling, be interpreted as elements of $I_{\epsilon_n}(S^4)$
(for an appropriate $\epsilon_n$) or $Al_1(S^4)$. We can also interprete them as elements of $Knot_{4,1}$ or $E_4$. Finally, we can construct
a hypersphere triangulated into flat simplices (instead of a
smooth hypersphere). It is easy to see that the number of
simplices will grow linearly with the length of the word $v_n$ in the
Baumslag-Gersten group that was used to construct, first, the balanced  
finite presentation of the trivial group and, then, a $4$-dimensional sphere.
Similarly, in the smooth case, ${vol^{1\over 4}\over inj}$, ${vol^{1\over 4}\over r}$ and $\vert K\vert diam^2$
will be bounded above by an exponential function of
$const$ $n$ for some $const$ (in fact, one can ensure much better
bounds, but we do not need this).

These hypersurfaces in $\mathbb{R}^5$ constructed using
the balanced presentations of the trivial group introduced
in [L] will be used in the proofs of all our results.
But note that in this construction one can alternatively
use another family of balanced presentations
of the trivial group with similar properties
that were independently discovered by Martin Bridson.
We were not aware of his work until this paper
had almost been finished. But after [L]
appeared on the arXiv Bridson e-mailed to us and wrote that he found
such finite presentations in 2003. Although they were mentioned
in his ICM-2006 talk ([B], p. 977), he has never publshed or posted any
details of his construction on the internet. 
Two weeks after the appearance of [L] his preprint [B2] has also appeared
on the arXiv.
Note that in his ICM-2006 talk Bridson expresses a hope
that such finite presentations of the trivial group
can be used to extend results of Nabutovsky
and Weinberger on the sublevel sets of diameter on moduli spaces
to dimension $4$ (which is something that we were not yet able to accomplish).
So, his work [B2] was also partially motivated by potential applications
that are similar in spirit to our results in this paper.
Even earlier, in the 90s, the second author attempted to prove the results
of this paper using balanced finite presentations of
the trivial group obtained from $B$ in the most obvious way, namely, by adding the second relation $w_n=t$. Yet he was not able
to verify that these balanced presentations have the desired properties.

{\bf 2.2. The filling length.} Following [N2]
we are going to use the following characteristic of simply-connected
closed Riemannian manifolds that measures how
``difficult" is to contract closed curves.
We define it as the supremum over all closed curves $\gamma$ of the ratio
${fl(\gamma)\over length(\gamma)}$, where the
filling length $fl(\gamma)$ denotes the infimum
over all homotopies $H=(\gamma_t)_{t\in [0,1]}$, $\gamma_0=\gamma$,
contracting $\gamma$ to a point (=constant curve) $\gamma_1$ of the maximal
length $\sup_t\ length(\gamma_t)$ of the closed curves arising during 
the homotopy $H$. We are going to denote this quantity by $Fl$ and regard
it as a functional on a considered space of (isometry classes) of
Riemannian metrics. 

To see that $Fl<+\infty$ first note that all sufficiently short curves
$\gamma$ can be contracted to a point without length increase (and,
therefore, $fl(\gamma)=length(\gamma)$). On the other hand for long curves
$\gamma$ we can choose any point $z$, and connecting $z$ with a sequence
of sufficiently close points on $\gamma$ by minimal geodesics reduce contraction
of $\gamma$ to consecutive contractions of triangles formed by a very short arc
of $\gamma$ and two minimal geodesics betwen $z$ and two very close
points on $\gamma$. Perimeters of these triangles are bounded by $2d+
\epsilon$, where $d$ is the diameter of the Riemannian manifold and $\epsilon$ is arbitrarily small. This easily implies that ${fl(\gamma)\over length(\gamma)}
\longrightarrow 1$ as $length(\gamma)\longrightarrow\infty$. This fact
was first noticed by M. Gromov ([Gr]). (The existence of the supremum for closed curves
of length $\leq 2d+\epsilon$ follows from the compactness of the set of 
Lipschitz curves of length $\leq 2d+\epsilon$ parametrized by the arclength.)
Gromov also introduced the term ``filling length" and the notation $fl$
(with a slightly different meaning than what we use here).

Note that $Fl$ can also be defined
for all simply-connected length spaces such that for some positive
$\epsilon$ all closed curves of length $\leq \epsilon$ can be contracted
to a point without length increase. So, in particular,
we can consider $Fl$ as a functional on the spaces of triangulations
of closed manifolds after we endow each simplex of the maximal dimension
by the metric of a regular flat simplex with side length $1$. (Actually, it is
easy to see that $Fl$ will not depend on the choice of side length here.)

Our observation is that if $S^4(v_n)$ is a (smooth or PL) sphere:
constructed starting from the word $v_n$ in the Baumslag-Gersten group
as expalined above (using either the idea from [L] or the idea from [B2])
then:

\begin{Pro}
The value of $Fl_n=Fl(S^4(v_n))$ grows faster than any finite tower of exponentials of $n$.
\end{Pro}

{\bf Proof.} Indeed, if not, then we can prove that the area of van Kampen diagrams
for generators of $S^4(v_n)$ will also be bounded by towers of exponentials
of $n$ of a fixed height. In the proof below we use the same notation $const$
for different constants that can be, in principle, evaluated.

The idea is that one can choose a way to represent each close curve $\gamma$
of length $\leq x$ by a word of length $\leq const\ n\ x$
so that if two closed
curves $\gamma_1$ and $\gamma_2$ are $const$-close, then
the corresponding words can be connected by a sequence of at most
$const\ x$ relations. In order to achieve this we first project $\gamma$
to the embedding of the representation complex of the balanced presentation
in $\mathbb{R}^5$. Recall, that $S^4(v_4)$ is the smoothed-out boundary
of a small tubular neighborhood of the representation complex, so this step
will increase the length by at most $const\ n$ factor. Denote the projection
of $\gamma$ to the embedded representation complex by $\tilde\gamma$.

Note that if $D$ is a Riemannian $2$-disc one can choose
a way to replace each arc with endpoints on $\partial D$ by a shortest
arc of $\partial D$ with the same endpoints. An ambiguity in the situation
when the distance between the endpoints along the boundary is equal to ${\vert\partial D\vert\over 2}$ leads to a discontinuity, and arcs corresponding
to two very close curves can together almost form the boundary of $D$.
Consider now one of the two $2$-cells in the representation complex
and a connected smooth arc $A$
in its interior with end points on its boundary.
The boundary of the $2$-cell has some self-intersections that appeared
as the result of taking the quotient map, when the cell was attached
to the $1$-skeleton. Yet we can canonically lift $A$ to the $2$-disc
with the same Riemannian metric in the interior and with nonself-intersecting
boundary, then extend $A$ to the boundary, replace it by a shortest
arc of the boundary with the same endpoints, and, finally, project this arc back
to the $1$-skeleton of the embedded representation complex.

Now take each component of the intersection of
$\tilde\gamma$ with $e_i$, $i=1,2$ and replace it by a minimal
arc in the boundary of $\partial e_i$ with the same endpoints as
explained above. This
will result in a length increase by at most $const\ n$ factor. 
The mentioned ambiguity in the case when the points of intersection
with $\partial e_i$ can be connected in the lift of $\partial e_i$
by two arcs of length ${\partial e_i\over 2}$ is
a source of discontinuity of this process.
If such a discontinuity arises for a pair of close curves $\gamma_1$, $\gamma_2$, then
the resulting arcs in the $1$-skeleton of the
representation complex will (almost) form the boundary of $e_i$.
Also, note that such ambiguity
(discontinuity) can arise only for a sufficiently long arc of $\gamma\bigcap e_i$, and, therefore, the number of such occurences for a
closed curve $\gamma$ is bounded by
$const\ x$. Now note that once $\gamma$ is replaced by a closed curve
in the $1$-skeleton of the representation complex, we can assign to it 
a word, and conclude that these words for sufficiently close curves can be transformed one into the other by an application of at most $const\ x$
relations.

Now we observe that a well-known and easy argument implies that
for each $\delta>0$ there
exists a $\delta$-net in the space of closed curves of length
$\leq x$ in the constructed Riemannian manifold of cardinality bounded
by $\exp(const\ {x\over\delta})$.
Therefore,
each closed curve can be contracted to a point by a discretized
homotopy that consists of at most $\exp(const\ x)$ ``jumps" of 
``length" $\leq const$ so that each ``jump" corresponds to a sequence of
not more than $const\ x$ applications of the relations for words
corresponding to the curves.

Finally, let $\gamma_0$ be a closed curve that represents one of the generators of
the considered finite presentation. It can be connected to a point
through curves of length $\leq Fl_n$ $length(\gamma_0)$. Therefore,
one obtains an at most exponential in $const\ Fl_n$
upper bound for the number of
the relations required to demonstrate that the generator, is, indeed, trivial.
This completes the proof of the proposition.

Let $M_0$ be any closed simply connected
$4$-dimensional Riemannian manifold.
We can form a Riemannian connected sum of $M_0$ with the spheres $S^4(v_n)$
in an obvious way and observe that $Fl$ for resulting Riemannian
manifolds grows faster than any tower of exponentials of $n$ of a fixed
height. 

{\bf 2.3. Proof of theorems.} It is now easy
to prove the main theorems for simply connected manifolds. In order to prove Theorem 1.4
recall that each bistellar transformation leads to a bi-Lipschitz homeomorphism
of the underlying simplicial complexes regarded as metric spaces, where each face of dimension four
is given the metric of the regular $4$-simplex with the side length $1$.
The Lipschitz constants for the map and its inverse do not exceed an absolute
constant $const$. Now note that in this situation $Fl$ cannot change by more than
the factor $const^2$. The value of $Fl$ for the boundary $\partial\Delta^5$
of the regular
$4$-simplex (endowed with the standard metric) is $1$. Therefore, the value
of $Fl$ on each triangulation of $S^4$ that can be connected with $\partial\Delta^5$ by at most $M$ bistellar transformation is at most $const^{2M}$.
This fact immediately implies the assertion of the theorem.

The proofs of the first three theorems are similar. The idea is to prove
that if the assertion does not hold, then $Fl_n$ is bounded above
by a tower of exponentials of $n$ of a fixed height, and this would
contradict the assertion of Proposition 2.1.

One can follow [N2] to finish the proof of Theorem 1.1. One starts from 
the observation that if two Riemannian structures in $I_\epsilon(M)$
are ${\epsilon\over 8.5}$-close (in the Gromov-Hausdorff metric), then
the values of $Fl$ can differ by a factor that does not exceed $1000000$
(Lemma 2 in [N2]). (The idea is that if $M_1$ and $M_2$ are close
Riemannian manifolds and one can contract any closed curve in $M_2$ through
not too long curves, one
can try to contract any closed curve $\gamma$ in $M_1$ by 1) discretizing
it, moving points to the closest points in $M_2$ and connecting them
by minimal geodesics, thereby obtaining a closed curve $\gamma_2$
that can be regarded as a ``transfer" of $\gamma$ to $M_2$; 2) Contracting
$\gamma_2$ through not too long curves in $M_2$; 3) Discretizing this
homotopy and transfering closed curves in the discretization back to
$M_1$; 4) Connecting the transfers of the 
nearest closed curves by homotopies in $M_1$, thus, obtaining a homotopy
contracting $\gamma$ in $M_1$.)

The second observation used in [N2] is that one can use the
well-known proof of the fact that $I_{\epsilon}(M)$ is precompact to
give an explicit uper bound of the form $\exp({const\over \epsilon^9})$
(in the four-dimensional case) for the cardinality of an $\epsilon/20$-net in $I_{\epsilon}(M)$. This estimate can then be used to conclude that any two Riemannian structures in the same connect component of $I_{\epsilon}(M)$ can be connected by a sequence
of $\epsilon/9$-long ``jumps", so that the number of jumps does not
exceed $\exp({const \over \epsilon^9})$ (see the proof of Lemma 3
in [N2]). Combining this estimate with the previous observation
we see that the ratio of values of $Fl$
at any two elements of the same connected
component of $I_{\epsilon}(M)$ 
is bounded by a double exponential of a power of
${1\over \epsilon}$ (and, thus, by a triple exponential
function of ${1\over\epsilon}$ for all sufficiently
small values of $\epsilon$). This estimate can be generalized to a stronger
equivalence relation on $I_{\epsilon}(M)$ than being in the same connected component,
namely, the transitive closure of the relation ``to be ${\epsilon\over 9}$-close in the Gromov-Hausdorff metric".

A comparison of these triply exponential upper bounds with lower bounds
for $Fl$ that grow faster than any tower of exponential of a fixed height
of $n$ yields the assertion of Theorem 1.1.

As it had been noticed, Theorem 1.2 would follow from the disconnectedness
of sublevel sets of the diameter $diam^{-1}((0,x])$ on $Al_1(M)$, and the
injectivity radius is bounded below by $\exp(-const\ x)$ on these sets.
Now one can use the same argument as in the proof of Theorem 1.1.

To prove Theorem 1.3 we can rescale the hypersurface
to have the value of the volume
equal to $1$. Now note that the definition of $\kappa$ implies that 
$\kappa\geq \vert k\vert$, where $k$ denotes any of the
principal curvatures of the hypersurface. 
This implies the obvious upper bound for the
absolute values of its sectional curvatures, when it is
regarded as a Riemannian manifold. It is not difficult to establish an
upper bound for the diameter of the hypersurface in the inner
metric (which immediately follows from Theorem 1.1 in [T]). Now
the Cheeger inequality implies an explicit lower bound for
the injectivity radius of the hypersurface that behaves as
$\exp(-const\ \kappa^{const})$, and we can prove the disconnectedness
of sublevel sets of $\kappa$ for an unbounded sequence of values
of $x$ exactly as we proved the disconnectedness of $I_{\epsilon}(M)$.

The proofs of Theorems 1.1, 1.2 and 1.4
in the case of a nonsimply connected manifold can be based
on the same ideas. We form a Riemannian connected sum of $M$ endowed
with some Riemannian metric with $S^4(v_n)$. Now $Fl$ is not defined,
but we can look at how much the length of the closed
curves corresponding to the
generators of the balanced presentations must be increased before they
can be contracted to a point. An argument similar to the proof of 
Proposition 2.1 implies that the growth of this quantity with $n$ is
faster than any tower of exponentials of a fixed height. On the other
hand the arguments in this section can be used to demonstrate
that the connectedness
assumptions imply a much slower growth of this quantity.

{\bf Acknowledgements.} This work has been partially supported from NSERC
Accelerator and Discovery Grants
of second named author (A.N.).

{\bf References.}

\noindent
[BN] V. Berestovskii, I. Nikolaev, ``Multidimensional generalized Riemannian
spaces", in ``Geometry iV", ed. Yu. G. Reshetnyak, Encyclopedia of
Mathematical Sciences, vol. 70, Springer, 1993.
\par\noindent
[BHP] W. Boone, W. Haken, and V. Poenaru, ``On Recursively Unsolvable Problems in Topology and Their Classification", Contributions to Mathematical Logic (H. Arnold Schmidt, K. Sch\"utte, and H. J. Thiele, eds.), North-Holland, Amsterdam, 1968.
\par\noindent
[B] M. Bridson, ``Non-positive curvature and complexity for finitely presented groups", Proceedings of the ICM, Madrid, Spain, 2006, 961-987,
European Mathematical Society, 2006.
\par\noindent
[B2] M. Bridson, ``The complexity of balanced presentations and the Andrews-Curtis
conjecture", arXiv:1504.04261
\par\noindent
[Ge] S. M. Gersten, Dehn functions and $l_1$-norms of finite presentations, Algorithms and Clas- sification in combinatorial group theory (Berkeley, CA, 1989) (G. Baumslag, C. Miller, eds.), Math. Sci. Res. Inst. Publ. 23, Springer-Verlag (1992),
195-224.
\par\noindent
[Gr] M. Gromov, ``Metric structures for Riemannian and non-Riemannian spaces", Birkhauser, 1998.
\par\noindent
[L] B. Lishak, ``Balanced finite presentations of the trivial group",
arXiv:1504.00418.
\par\noindent
[N1] A. Nabutovsky, ``Non-recursive functions, knots ``with thick ropes",
and self-clenching ``thick" hypersurfaces", Comm. Pure Appl. Math. 48(1995), 381-428.
\par\noindent
[N2] A. Nabutovsky, ``Disconnectedness of sublevel sets of some Riemannian functionals", Gem. Funct. Analysis (GAFA), 6(1996), 703-725.
\par\noindent
[N3] A. Nabutovsky, ``Geometry of the space of triangulations of a compact manifold", Comm. Math. Phys. 181(1996), 303-330.
\par\noindent
[N4] A. Nabutovsky,  ``Combinatorics of the space of Riemannian structures and
logic phenomena of Euclidean Quantum Gravity", in ``Perspectives in Riemannian
Geometry", ed. by V. Apostolov et al., CRM Proceedings and Lecture Notes,
vol. 40, 223-248, AMS, Providence, RI, 2006.
\par\noindent
[N5] A. Nabutovsky,  ``Morse landscapes of Riemannian functionals and related
problems", Proceedings of ICM-2010, vol. 2, 862-881, Hindustan Book Agency,
New Delhi, 2010.
\par\noindent
[NW1] A. Nabutovsky, S. Weinberger, ``Variational problems for Riemannian functionals and arithmetic groups", Publications d'IHES, 92(2000), 5-62.
\par\noindent
[NW2] A. Nabutovsky, S. Weinberger, ``The fractal nature of Riem/Diff I",
Geom. Dedicata 101(2003), 1-54.
\par\noindent
[Ni] I. Nikolaev, ``Bounded curvature closure of the set of compact
Riemannian manifolds", Bull. of the AMS 24(1992), 171-177.
\par\noindent
[M] J. Milnor, ``Lectures on the $h$-cobordism theorem", Princeton University Press, Princeton, NJ, 1965.
\par\noindent
[P] U. Pachner, ``P.L. homeomorphic manifolds are equivalent by elementary
shellings", Europ. J. Combinatorics 12(1991), 129-145.
\par\noindent
[Pl] A. N. Platonov An isoperimetric function of the Baumslag-Gersten group, Moscow Univ. Math.
Bull. 59 (2004), no 3, 12-17.
\par\noindent
[Sh 1] M. A. Stan'ko ``On the Markov theorem on algorithmic
nonrecognizability of manifolds", J. of Mathematical Sci., 146(2007), 5622-5623.
\par\noindent
[Sh 2] M. A. Stan'ko ``The Markov theorem and algorithmically
unrecognizable combinatorial manifolds", Izv. Akad. Nauk SSSR, Ser. Math., 68(2004), 207-224.
\par\noindent
[T] P. Topping,  ``Relating diameter and mean curvature for submanifolds of
Euclidean space", Comment. Math. Helv. 83(2008), 539-546.
\par\noindent
[W] S. Weinberger, ``Computers, rigidity and moduli", Princeton University Press, 2004.

\end{document}